\documentclass[a4paper,11pt]{article}
\usepackage{amsmath, amssymb, amscd, amsthm, euscript, color}

\newtheorem{thm}{Theorem}[section]
\newtheorem{cor}[thm]{Corollary}

\newtheorem{lem}[thm]{Lemma}

\newtheorem{rem}[thm]{Remark}

\title{Generating  the M\"obius group with involution conjugacy classes}
\author{Ara Basmajian\footnote{Supported in part by PSC-CUNY Grant 627 14-00 40.} and Karan Puri\footnote{The author is indebted to his advisor, Jane Gilman and would like to thank Rutgers-Newark for its support. It is a pleasure to thank Mark Feighn, Lee Mosher and Jacob Sturm for helpful conversations and encouragement.}}

\begin{document}
\bibliographystyle{plain}

\maketitle

\begin{abstract}
 A  {\it $k$-involution}  is an involution with a fixed point set of codimension $k$. The conjugacy class of such an involution, denoted $S_k$,  generates  $\text{M\"ob}(n)$-the the group of isometries of hyperbolic $n$-space-if $k$ is odd, and its orientation preserving  subgroup if $k$ is even.  In this paper, we supply effective lower and  upper bounds for
the  $S_k$ word  length of    $\text{M\"ob}(n)$ if $k$ is  odd, and the
$S_k$ word  length of   $\text{M\"ob}^+(n)$,   if $k$ is  even.  As a consequence,  for a fixed codimension  $k$   the length of
$\text{M\"ob}^{+}(n)$ with respect to $S_k$, $k$ even,  grows linearly with $n$ with the same statement  holding in the  odd case.    Moreover,  the percentage of  involution  conjugacy classes  for which  $\text{M\"ob}^{+}(n)$ has length two approaches zero, as $n$ approaches infinity.
\end{abstract}

\section{Introduction and results.}

  Let $G$ be a group and $S$ a set of symmetric  generators for a supergroup of $G$; $S$ is not necessarily a subset of $G$ but every element can be written as a product of elements from $S$.  For $g \in G$, the {\it length} of $g$ with respect to $S$ (or $S$-length)  is the minimal number of elements of $S$ needed to express $g$ as their product.  The supremum  over all  group  element lengths is called the {\it length  of $G$  with respect to $S$}  (or simply the $S$-length of $G$), and is denoted $|G|$.

We are interested in the  set $S_k \subset  \text{M\"ob}(n)$ of  involutions  with a codimension $k$ fixed point set acting on  hyperbolic space,  $\mathbb{H}^n$.

\newpage

\begin{thm}\label{thm:bounds}    Let  $n \geq 2$ and   $k=1,2,..., n-1$.

\begin{itemize}
  \item    If $k$ is even, $S_k$ generates   $\text{M\"ob}^{+}(n)$ and satisfies,

\begin{equation}
 \frac{n(n+1)}{2k(n-k+1)} \leq  |\text{M\"ob}^{+}(n)|_k  \leq 2n+4. \label{eq:bounds}
\end{equation}

  \item    If $k$ is odd, $S_k$ generates   $\text{M\"ob}(n)$ and satisfies,

\begin{equation}
 \frac{n(n+1)}{2k(n-k+1)} \leq  |\text{M\"ob}^{+}(n)|_k   \leq 2n+2+k. \label{eq:bounds}
\end{equation}
\end{itemize}
\end{thm}
where $| \cdot |_k$  denotes $S_k$-length.

In particular  we have,

\begin{cor}[Linear growth]
 For a fixed codimension $k$,

 \begin{itemize}

  \item   If $k$ is even, $|\text{M\"ob}^{+}(n)|_k \asymp n$.

   \item  If $k$ is odd,  $|\text{M\"ob}(n)|_k \asymp n$.

\end{itemize}

\end{cor}

 Our main tool is  a dimension count  which yields the fact that
the set of elements of the form $\alpha_1...\alpha_m$, where $\alpha_i$  is a $k$-involution
and  $m<\frac{n(n+1)}{2k(n-k+1)}$,  has measure zero in    $ \text{M\"ob}(n)$ (see  theorem \ref{thm:measurezero}  and corollary \ref{cor:Mobgen}).

We next consider the percentage of involution  conjugacy classes  for which
$|\text{M\"ob}^+(n)|=2$.  More precisely, define

\begin{equation}
\Phi(n)=\frac{|\{k: \text{the $S_k$-length of $\text{M\"ob}^+(n)$ is 2}\}|}{|\{\text{involution conjugacy classes in dimension $n$}\}|}.
\end{equation}

 Recall that the conjugacy class  of  an involution is determined  by the dimension of its  (totally geodesic) fixed point set.  Thus a $k$-involution determines a conjugacy class,  and the denominator above is  $n$.
   Setting  $|\text{M\"ob}^+(n)|=2$   and using  the lower bound in    inequality  (\ref{eq:bounds}), after a     straightforward computation  to find bounds on $k$ in terms of $n$ we have,

 \begin{cor}
$ \Phi(n)= O(n^{-\frac{1}{2}})$.
  \end{cor}

In the paper    \cite{B-M} the authors show that in each dimension  $n$  there exists an involution conjugacy class $S_k$  for which  $|\text{M\"ob}^+(n)|=2$. When $n$ is even $k$ may be taken to be $\frac{n}{2}$, and when $n$ is odd,  $k$ can be taken to be  $\frac{n+1}{2}$.
 The results  of this paper show  that   away from the  middle codimensions  (relative to $n$)
one cannot expect the length of  $\text{M\"ob}^+(n)$  to be small. Factoring of isometries from a geometric  viewpoint  in hyperbolic 4-space  is also studied in \cite{P} where a slightly different formulation and generalization of two and three dimensional half-turns is given.

    Hyperbolic  $n$-space is denoted $\mathbb{H}^n$.
    The orientation preserving  isometries of  $\mathbb{H}^n$   (the M\"obius group)  is   $\text{M\"ob}^{+}(n)$,  and the full group is    $\text{M\"ob}(n)$.  An involution is an order two isometry of    $\mathbb{H}^n$ and a  {\it $k$-involution} is an isometry with a fixed point set of codimension $k$. A {\it reflection} is an involution with a codimension one  fixed point set, and a {\it half-turn} is an involution with a codimension two fixed point set.     The involution  is  orientation reversing if and only if  the codimension of the fixed point set  is  odd. For  $k=1,2,3,...,n$, let $S_k$ be the set  (conjugacy class) of  $k$- involutions.  For the basics on hyperbolic space and  its isometry group we refer to Maskit or Ratcliffe (\cite{M},\cite{R}). A good background for $n$-dimensional hyperbolic geometry and some lower dimensional factorization results are given in \cite{Wi}. For standard material on  differential topology and   Lie groups  the reader is referred to \cite{S}  and  \cite{W}.

 We will use the well-known  facts that  the dimension of $O(n)$, as well as $SO(n)$, is  $\frac{n(n-1)}{2}$ and the dimension of  $\text{M\"ob}^{+}(n)$, as well as   $\text{M\"ob}(n)$,  is $\frac{n(n+1)}{2}$.

 The paper is organized as follows.  Section   (\ref{sec:$k$-involutions}) contains the proofs that
conjugacy classes of  involutions  generate  the M\"obius group   as well as upper bounds on word length.   In section  (\ref{sec:Involutions and the space of totally geodesic subspaces}),  we show that the space of $k$-involutions  is a submanifold of  $\text{M\"ob}(n)$ having dimension $k(n-k+1)$. Finally, we prove   theorem  \ref{thm:bounds}  in section (\ref{sec:Bounds for the $S_k$-length}).


 \section{$k$-involutions in the orthogonal  and  M\"obius groups} \label{sec:$k$-involutions}

Throughout this section,  we   fix  an integer $n \geq 2$ and an integer  $k=1,2,..., n-1$.
 The case $k=n$ is  excluded since $S_n  \cap O(n)$  has only one element  and hence does not  generate  the orthogonal group.

  This section is devoted to proving,

 \begin{thm}
 \label{theorem:orthogonalgen}   $S_k \cap  O(n)$ generates  $O(n)$ if $k$ is odd, and generates $SO(n)$  when $k$ is even.
 Furthermore for $g \in O(n)$,

 \begin{equation}  |g| \leq \left\{
  \begin{array}{l}
     2n,   \qquad \qquad   \text{if  $g$ is orientation preserving}   \\
     2n-2+k,  \,\,  \text{if $g$ is orientation  reversing} \\
\end{array}
\right\},
\end{equation}
where  $|g|$ is the   $S_k \cap  O(n)$-length of $g$. \end{thm}

 \begin{cor}
 \label{cor:Mobgen}   $S_k$ generates  $\text{M\"ob}(n)$   if $k$ is odd,  and generates   $\text{M\"ob}^{+}(n)$ when $k$ is even.  Furthermore for $g \in  \text{M\"ob}(n)$,

 \begin{equation}  |g| \leq \left\{
  \begin{array}{l}
     2n+4,   \qquad \qquad   \text{if  $g$ is orientation preserving}   \\
     2n+2+k,  \;\;\;\;\;\;\;\,  \text{if $g$ is orientation  reversing} \\
\end{array}
\right\},
\end{equation}
where  $|g|$ is the   $S_k$-length of $g$.   \end{cor}

\begin{rem}

In both theorem   \ref{theorem:orthogonalgen} and corollary   \ref{cor:Mobgen}, we note  that  $k$ is necessarily odd when $g$ is orientation reversing.

  \end{rem}

The  stabilizer of any point in $\mathbb{H}^n$ has a natural  identification with $O(n)$. We fix such a copy of  $O(n) \subset  \text{M\"ob}(n)$.

 Denote the $n\times n$ diagonal matrices with $k$ entries being $-1$ and $n-k$ being $1$  by
 $ \mathcal{D}(n,k)$.  Since  an involution in $O(n)$  is  $O(n)-$conjugate to  a diagonal matrix, it is immediate  that a  $k$-involution in $O(n)$ is conjugate  to a diagonal  matrix  in
   $ \mathcal{D}(n,k)$. There are    $ \left(\begin{array}{c}n \\k\end{array}\right)$   such matrices.

\begin{lem} \label{lemma:products}

Assume $n\ge 2$ and $k=1,...,n-1$.

\begin{enumerate}

\item If $k$ is odd, then any element of
 $ \mathcal{D}(n,1)$   can be written as the product of  $k$ elements of $\mathcal{D}(n,k)$.

\item
Any element of $\mathcal{D}(n,2)$ can be written as the product of two elements of $\mathcal{D}(n,k)$.

\end{enumerate}
\end{lem}

\begin{proof}  For ease of notation, we identify the diagonal matrices of size $n$  having  $\pm 1$  entries with the group $\mathbb{Z}_2^n$. That is,
    $\bigcup_{k=0}^{n} \mathcal{D}(n,k)= \mathbb{Z}_2^n$. We write an element of   $\mathbb{Z}_2^n$ as a vector with the obvious component-wise multiplication in  $\mathbb{Z}_2$.

 To prove item (1),  Consider $A$=$[-1,1,1,...,1] \in\mathcal{D}(n,1)$. It suffices to show that $A$ can be written as the desired product. For $i=1,..., k$, let $C_i\in\mathcal{D}(n,k)$  with   $j$-th
 component  being,

 \begin{equation} C_i^j=\left\{
  \begin{array}{l}
     1,   \,\,  \text{if  $j=i+1$ or $ k+2\leq j \leq n$}   \\
     -1,  \,\,  \text{if $1 \leq j \leq  k+1$ and $j \neq  i+1$} \\
\end{array}
\right\}
\end{equation}
  Then  $A=\displaystyle\prod_{i=1}^{k} C_i$  and we have the desired decomposition of $A$.

To prove item (2),   consider $A=[-1,-1,1,...,1]\in\mathcal{D}(n,2)$. It suffices to show that $A$ can be written as the desired product. Let $R\in\mathcal{D}(n,k)$ be  such that  its $j$-th component is,

 \begin{equation}  R^j=\left\{
  \begin{array}{l}
     1,   \,\,  \text{if  $j=1$ or $ k+2 \leq j \leq n$}   \\
     -1,  \,\,  \text{if $2 \leq j \leq k+1$} \\
\end{array}
\right\}
\end{equation}

and let $S \in\mathcal{D}(n,k)$ have $j$th entry

 \begin{equation}  S^j=\left\{
  \begin{array}{l}
     1,   \,\,  \text{if  $j=2$ or $ k+2 \leq j \leq n$}   \\
     -1,  \,\,  \text{otherwise}. \\
\end{array}
\right\}
\end{equation}

 Then  $RS=A$ and we are finished with the proof of item (2).

  \end{proof}

  \begin{lem}

  \label{lemma:generalcodim2}  Let $a$ and $b$ be reflections in hyperplanes $\alpha$ and $\beta$ in $\mathbb{H}^n (n\geq 3)$ and let $g=ab$. Then, there exist half-turns $h$ and $k$ such that $g=hk$.

  \end{lem}

  \begin{proof}

  Consider the upper half-space model of $\mathbb{H}^n$. Let $\alpha\cap\hat{\mathbb{R}}^{n-1}=\tilde{\alpha}$ and $\beta\cap\hat{\mathbb{R}}^{n-1}=\tilde{\beta}$. Then $\tilde{\alpha}$ and $\tilde{\beta}$ are $(n-2)$-spheres in $\hat{\mathbb{R}}^{n-1}$. We may assume that neither $\tilde{\alpha}$ nor $\tilde{\beta}$ contains the point at infinity ($\infty$). Consider the unique circle $\rho$, through $\infty$ and each of the centers of $\tilde{\alpha}$ and $\tilde{\beta}$. It is clear that any $(n-2)$-sphere  containing  $\rho$ is orthogonal to each of $\tilde{\alpha}$ and $\tilde{\beta}$. Let $\tilde{\gamma}$ be one such $(n-2)$-sphere.

  Then, $\tilde{\gamma}=\gamma\cap\hat{\mathbb{R}}^{n-1}$, where $\gamma$ is a hyperplane in $\mathbb{H}^n$ which is orthogonal to each of $\alpha$ and $\beta$. Let $c$ denote reflection in $\gamma$. Then, $h=ac$ and $k=cb$ are half-turns in $\mathbb{H}^n$ such that $hk=(ac)(cb)=accb=ab=g$.

  \end{proof}

 \begin{proof}[Proof of theorem  \ref{theorem:orthogonalgen}]  Fix $k=1,...,n-1$.  Using the block diagonal form an element in $g \in SO(n)$,
it is easy to see  that  an element    $g \in SO(n)$  can be written as  a product  $\rho_1... \rho_m$,  where  $\rho_i \in  S_1 \cap O(n)$, $m$ is even,  and $m$ is  at most $n$.   Now,   using    lemma  \ref{lemma:generalcodim2}
 we  write $g$ as a product of $m$ half-turns.   Of course, the half-turns are     $O(n)$-conjugate
 to  a diagonal matrix  in   $ \mathcal{D}(n,2)$ and hence using  item (2) of lemma
 \ref{lemma:products}, we can  write  $g$ as the product of at most  $2n$ elements in
  $S_k  \cap O(n)$.

 If $g \in  O(n)-SO(n)$, then  $g=\rho_1... \rho_m$, where $m$ is odd and at most $n$.  Note that  it must be that  $k$ is  odd.   As above we write $\rho_1... \rho_{m-1}$ as the  product of  at most   $2n-2$ elements  in $S_k \cap O(n)$.  The reflection $\rho_m$, using  item (1) of lemma     \ref{lemma:products}, can be written as the product of $k$  elements  in   $S_k \cap O(n)$. Thus for such an element  $g$, $|g| \leq  2n-2+k$.
   \end{proof}

\begin{proof} [Proof of  corollary   \ref{cor:Mobgen}]\
For  $g \in  \text{M\"{o}b}(n)$, it is well-known  that $g=\Phi  \tau\sigma $, where $\sigma$ and
 $\tau$ are reflections,  and $\Phi$ is  an  element of $O(n)$. Moreover $g$ is orientation preserving if and only if $\Phi \in    SO(n)$.  Using  lemma  \ref{lemma:generalcodim2},  we can replace  $\tau\sigma$ by the product of two half-turns which  by lemma
 \ref{lemma:products}  can be written as the product of 4 elements in  $S_k$. The  corollary now follows from theorem   \ref{theorem:orthogonalgen}.
\end{proof}


\section{Involutions and the space of totally geodesic subspaces of $\mathbb{H}^n$. }
\label{sec:Involutions and the space of totally geodesic subspaces}

Throughout this section,  we   fix  an integer $n \geq 2$ and an integer  $k=1,2,..., n-1$.

 \begin{lem}\label{lem:diffstructureonplanes}   $S_k \subset \text{M\"ob}(n)$ is a  (connected) differentiable submanifold  of dimension $k(n-k+1)$.  \end{lem}

 \begin{proof}   Set $G=\text{M\"ob}(n)$. Fix $\alpha \in S_k \subset  \text{M\"ob}(n)$ and  denote its fixed point set by $\pi$, an $ (n-k)$ dimensional plane. Consider the  smooth conjugation action of the Lie group $G$ on itself. Namely,  $g \cdot  f =gfg^{-1}$. Since an orbit of a Lie group action is a submanifold, we have that  the $G$-orbit of $\alpha$, that is $S_k$,  is a submanifold of  $\text{M\"ob}(n)$. Furthermore, the map  from $G$ to $G$ given by  $g \mapsto g\alpha g^{-1}$, induces a one-to-one smooth map  from $G/K$ onto  $S_k$, where $K=\text{Stab}_G(\alpha)$.
 (Note that $K$ is a closed subgroup of $G$). Next observe  that  $\text{Stab}_G(\alpha)=\text{Stab}_G(\pi)$ and consider the map,

 \begin{equation}
\Phi :  \text{Stab}_G(\pi) \rightarrow  \text{M\"ob}(n-k)
 \end{equation}
 given by $g \mapsto  g|_{\pi}$. This is a surjective map with kernel being isomorphic to
  $O(k) \leq  \text{Stab}_G(\pi) $. Hence,  $\text{Stab}_G(\pi) /O(k)$ is isomorphic to  $\text{M\"ob}(n-k)$
  and thus,

   \begin{equation}
\text{dim}(K)=\text{dim}(\text{Stab}_G(\alpha))= \text{dim}(\text{M\"ob}(n-k))+ \text{dim}(O(k))
  \end{equation}

  Thus we have

 \begin{equation}
 \text{dim}(S_k)=\text{dim}(G)-\text{dim}(K)
 =\text{dim}(G)-\text{dim}(\text{M\"ob}(n-k))-\text{dim}(O(k))
  \end{equation}

Now plugging in the various quantities and simplifying yields the dimension of $S_k$ to   be
$k(n-k+1)$.

  \end{proof}

  For $k=1,..,n-1$, let  $\mathcal{G}_k$ denote the space of  $k$-planes  (that is, $k$-dimensional totally geodesic subspaces)
   in  $\mathbb{H}^n$.    The boundary (at infinity) of a $k$-plane is a round  $(k-1)$ sphere. The space of
   $(k-1)$ spheres with the Gromov-Hausdorff topology induces a natural  topology on
   $\mathcal{G}_k$.

  \begin{cor}
 $\mathcal{G}_k$ is a  differentiable manifold  of dimension $(n-k)(k+1)$.
  \end{cor}

  \begin{proof}
  Consider the map, $ \mathcal{G}_k \rightarrow \mathcal{S}_{n-k}$ given by taking  the  $k$-plane $\pi$  to the   $(n-k)$-involution  with fixed point set $\pi$. This map is a homeomorphism (needs to be checked).
  Pulling back the differentiable structure from $\mathbb{S}_{n-k}$,   $\mathcal{G}_k $ becomes a differentiable manifold whose dimension by lemma   (\ref{lem:diffstructureonplanes})  is  $(n-k)(k+1)$.

  \end{proof}


\section{Bounds for the $S_k$-length of the M\"obius group}\label{sec:Bounds for the $S_k$-length}

Given  a subset $J \subseteq  \{1,...,n-1\}$, let $S$ be
the generating set  $S=\cup_{k\in J} S_k$ and set $M=M(S)=\text{max}_{c\in J}\{\text{dim}(S_{k})\}=
\text{max}_{k\in J}\{ k(n-k+1)\}$.

\begin{thm}\label{thm:measurezero}
 Except for  a set of measure zero no element of     $\text{M\"ob}(n)$ can be written as a product
 $\alpha_1...\alpha_m$, where $\alpha_i \in S$,  and  $m<\frac{n(n+1)}{2M}$.
 \end{thm}

\begin{proof}  Consider the manifold which is the  $m$-fold product of the M\"obius group.
 Given  a sequence  $\{k_1,...,k_m\}$ of $m$-elements from $J$ (repetition is allowed),  consider the mapping  $\Psi :S_{k_1}\times...\times S_{k_m} \rightarrow    \text{M\"ob}(n)$  which assigns   the $m$-tuple of  ordered  $k_i$-involutions $(\alpha_1, ... ,\alpha_m)$    to  the product $\alpha_1...\alpha_m$.   This is a smooth mapping between manifolds.  The dimension of   $S_{k_1}\times...\times S_{k_m} $  is   bounded from above  by $mM$  which is by assumption  less than the dimension of   $\text{M\"ob}(n)=\frac{n(n+1)}{2}$.    Hence $\Psi$ is a smooth mapping from a  manifold of lower dimension to one of higher dimension.  It is a standard fact   that the image of  a  smooth map from  a  manifold  of lower dimension to one of higher dimension   has measure zero.

Finally, there are a finite number of  (namely, $\sum_{i=1}^{[\frac{n(n+1)}{2M}]} |J|^i$)  sequences from  $J$  of length less than or equal to $\frac{n(n+1)}{2M}$, and hence a finite number of maps $\Psi$ above.  Thus the set of elements in   $\text{M\"ob}(n)$ that  are in the images  of such maps $\Psi$ is the finite union of sets of  measure zero, hence has measure zero. This precisely says  that the elements of  $\text{M\"ob}(n)$ that  can be written as a product of  at most $m$ elements from $S$ has measure zero.

\end{proof}

{\bf Theorem 1.1.}   Let  $n \geq 2$ and   $k=1,2,..., n-1$.

\begin{itemize}
  \item    If $k$ is even, $S_k$ generates   $\text{M\"ob}^{+}(n)$ and satisfies,

\begin{equation}
 \frac{n(n+1)}{2k(n-k+1)} \leq  |\text{M\"ob}^{+}(n)|_k  \leq 2n+4. \label{eq:bounds}
\end{equation}

  \item    If $k$ is odd , $S_k$ generates   $\text{M\"ob}(n)$ and satisfies,

\begin{equation}
 \frac{n(n+1)}{2k(n-k+1)} \leq  |\text{M\"ob}^{+}(n)|_k   \leq 2n+2+k. \label{eq:bounds}
\end{equation}
\end{itemize}

\begin{proof}   The upper bound in either the even or odd case  follows from corollary \ref{cor:Mobgen}.  For the lower bound,
 if all $k_i=k$ then  $M=k(n-k+1)$.  Now if  $ |g|  < \frac{n(n+1)}{2k(n-k+1)} $,      for all  $g \in \text{M\"ob}(n)$ then theorem    (\ref{thm:measurezero}) is contradicted.
\end{proof}


\author{Ara Basmajian, Department of Mathematics, Graduate Center and Hunter College, CUNY, abasmajian@gc.cuny.edu}

\bigskip

\author{Karan Puri, Department of Mathematics, Queensborough Community College, CUNY, kpuri@qcc.cuny.edu}


\begin{thebibliography}{999}
\bibitem{B-M} A. Basmajian  and B. Maskit \textit{Space form isometries as commutators and products of  involutions}, preprint.

\bibitem{M}  B. Maskit,  \textit{Kleinian groups}. Grundlehren der Mathematischen Wissenschaften [Fundamental Principles of Mathematical Sciences], 287. Springer-Verlag, Berlin, 1988. xiv+326 pp.




\bibitem{P} K.M. Puri, \textit{Factorization of isometries of hyperbolic 4-space and a discreteness
condition},  Thesis, Rutgers University, 2009.

\bibitem{R}  J. G.  Ratcliffe,  \textit{Foundations of hyperbolic manifolds}. Second edition. Graduate Texts in Mathematics, 149. Springer, New York, 2006. xii+779 pp.

\bibitem{S}   M.  Spivak,  \textit{A comprehensive introduction to differential geometry}. Vol. I. Second edition. Publish or Perish, Inc., Wilmington, Del., 1979. xiv+668 pp.

\bibitem{Wi} J. B. Wilker, \textit{Inversive Geometry}. The Geometric Vein, 379--442, Springer, New York-Berlin 1981.

\bibitem{W} J. Wolf,    \textit{Spaces of Constant Curvature},  Publish or Perish,  5 ed.

\end{thebibliography}
\end{document}